%svt.tex: On the Number of Singular Vector Tuples of Hyper-Cubical Tensors
%by Shalosh B. Ekhad and Doron Zeilberger
%Plain TeX
%begin macros

\baselineskip=14pt
\parskip=10pt

\magnification=\magstephalf

\def\1{{\overline{1}}}
\def\2{{\overline{2}}}
\parindent=0pt
\overfullrule=0in

\def\frac#1#2{{#1 \over #2}}
%\headline={\rm  \ifodd\pageno  \RightHead  \else  \LeftHead  \fi}
%\def\RightHead{\centerline{
%Title
%}}
%\def\LeftHead{ \centerline{Doron Zeilberger}}
%end macros

\bf
\centerline
{
On the Number of Singular Vector Tuples of Hyper-Cubical Tensors
}
\rm
\bigskip
\centerline
{\it By Shalosh B. EKHAD and Doron ZEILBERGER}
\bigskip
\qquad 

A week ago, Bernd Sturmfels [St] gave a fascinating Colloquium talk, here at Rutgers,
where, among many other interesting facts, he mentioned the
following theorem of Shmuel Friedland and Giorgio Ottaviani ([FO]).

{\bf Theorem ([F0], Theorem 1)}. The number of simple singular vector tuples of a generic \hfill\break
$m_1 \times \cdots \times m_d$ ($d$-dimensional) tensor
equals the coefficient of  $\prod_{i=1}^{d} t_i^{m_i-1}$ in the polynomial
$$
\prod_{i=1}^{d} \frac{ {\hat{t_i}}^{m_i}-t_i^{m_i} }{\hat{t_i}-t_i} \quad , \quad  \hat{t_i}=\left ( \sum_{j=1}^{d} t_j  \right ) - t_i \quad .
$$
Let's call this number $c(m_1, \dots, m_d)$.

We first observe that the generating function of this $d$-dimensional multi-sequence is 
a nice rational function.

{\bf Proposition 1.} Let $e_i(x_1, \dots, x_d)$ be the elementary symmetric function of the indeterminates $x_1, \dots, x_d$, 
of degree $i$. We have:
$$
\sum_{m_1=0}^{\infty} \dots  \sum_{m_d=0}^{\infty}
c(m_1, \dots, m_d) x_1^{m_1} \dots x_d^{m_d}=
\prod_{i=1}^{d} x_i \left ( \prod_{i=1}^{d} (1-x_i) \right )^{-1} \left ( 1-\sum_{i=2}^{d} (i-1) e_i(x_1, \dots , x_d) \right )^{-1} \quad .
$$

{\bf Proof}: Since
$$
\frac{ {\hat{t_i}}^{m_i}-t_i^{m_i} }{\hat{t_i}-t_i} =\sum_{k_i=0}^{m_i-1} {\hat{t_i}}^{k_i} t_i^{m_i-1-k_i} \quad ,
$$
$c(m_1, \dots, m_k)$ is the coefficient of  $\prod_{i=1}^{d} t_i^{m_i-1}$ in
$$
\sum_{k_1=0}^{m_1-1} \dots  \sum_{k_d=0}^{m_d-1}
\,\, \prod_{i=1}^{d} {\hat{t_i}}^{k_i} t_i^{m_i-1-k_i} \quad .
$$
Hence
$$
c(m_1, \dots , m_d) \, = \,
\sum_{k_1=0}^{m_1-1} \dots  \sum_{k_d=0}^{m_d-1}
ConstantTermOf \prod_{i=1}^{d} \, {\hat{t_i}}^{k_i} t_i^{-k_i} \quad .
$$
Let
$$
f(k_1, \dots, k_d):=ConstantTermOf \prod_{i=1}^{d} \, {\hat{t_i}}^{k_i} t_i^{-k_i}  
\, = \,
CoeffOf \prod_{i=1}^d {t_i}^{k_i} \quad in \quad  \prod_{i=1}^d \left ( \sum_{j=1}^{i-1} t_j \, + \,  \sum_{j=i+1}^{d} t_j \right )^{k_i} \quad .
$$
By the celebrated {\bf MacMahon Master Theorem} ([M], Section III, Chapter II, p. 93ff) 
(with the $d \times d$ matrix that is all $1$'s except $0$ in the diagonal), we have
$$
\sum_{k_1=0}^{\infty} \dots  \sum_{k_d=0}^{\infty}
f(k_1, \dots, k_d) x_1^{k_1} \dots x_d^{k_d}=
\left ( 1-\sum_{i=2}^{d} (i-1) e_i(x_1, \dots , x_d) \right )^{-1} \quad .
$$
Since
$$
c(m_1, \dots , m_d) \, = \, \sum_{k_1=0}^{m_1-1} \dots  \sum_{k_d=0}^{m_d-1} f(k_1, \dots , k_d) \quad ,
$$
the proposition follows by straightforward generatingfunctionology.

The fact that the generating function of $c(m_1, \dots , m_d)$  is a rational function
is equivalent to it satisfying a certain partial linear recurrence with constant coefficients, 
easily deduced from the generating function.
Combined with the fact that both $c(m_1, \dots, m_d)$ and $f(k_1, \dots, k_d)$ are symmetric,
enabled us to efficiently compute many values. It also follows 
(for example using Wilf-Zeilberger algorithmic proof theory, efficiently implemented in [AZ])
that the diagonal sequences
$$
C_d(n):=c(n, \dots , n)  \quad \quad [ \, n \quad repeated \quad d \quad times \,]  \quad ,
$$
are {\bf holonomic}, alias {\bf P-recursive}, that means that for each $d$, the sequence
$C_d(n)$ satisfies {\it some} homogeneous linear recurrence with polynomial coefficients. 
While one can use the method of [AZ], it is more efficient, since we know {\it a priori} that
such a recurrence exists, to generate sufficiently many terms and then
{\bf guess} the recurrence. Using this method, we got the following proposition.

{\bf Proposition 2.} The sequence $C_3(n)=c(n,n,n)$ satisfies the following fifth-order
linear recurrence equation with polynomial coefficients.
$$
2\, \left( n+2 \right)  \left( 245\,{n}^{4}+3094\,{n}^
{3}+14447\,{n}^{2}+29474\,n+22100 \right)  \left( n+1
 \right) ^{2}C_3 \left( n \right) 
$$
$$
- \left( n+2 \right) 
 \left( 21805\,{n}^{6}+330981\,{n}^{5}+2012733\,{n}^{4
}+6230951\,{n}^{3}+10263446\,{n}^{2}+8425060\,n+
2639760 \right) C_3 \left( n+1 \right) 
$$
$$
+ \left( -13230\,
{n}^{7}-249641\,{n}^{6}-1998705\,{n}^{5}-8785333\,{n}^
{4}-22847777\,{n}^{3}-35069178\,{n}^{2}-29331496\,n-
10279296 \right) \cdot
$$
$$
C_3 \left( n+2 \right) 
$$
$$
+ \left( 21560\,
{n}^{7}+413637\,{n}^{6}+3343917\,{n}^{5}+14735333\,{n}
^{4}+38132651\,{n}^{3}+57777574\,{n}^{2}+47273504\,n+
16026528 \right) \cdot
$$
$$
C_3 \left( n+3 \right) 
$$
$$
- \left( n+4
 \right)  \left( 4410\,{n}^{6}+70147\,{n}^{5}+452903\,
{n}^{4}+1516515\,{n}^{3}+2769127\,{n}^{2}+2601986\,n+
975888 \right) C_3 \left( n+4 \right) 
$$
$$
+ \left( n+5
 \right)  \left( n+4 \right)  \left( n+3 \right) 
 \left( 245\,{n}^{4}+2114\,{n}^{3}+6635\,{n}^{2}+8882
\,n+4224 \right) C_3 \left( n+5 \right) =0 \quad ,
$$
subject to the initial conditions
$$
C_3(1) = 1 \quad , \quad C_3(2) = 6 \quad , \quad C_3(3) = 37 \quad , \quad C_3(4) = 240 \quad , \quad C_3(5) = 1621 \quad .
$$

Using the methods of [WZ] and [Z], we found the following asymptotic formula.

{\bf Proposition 3.}
$$
C_3(n) \, \sim \, \frac{2}{\sqrt{3} \, \pi} \, 8^n  \, \cdot \, n^{-1} \cdot
$$
$$
\left( 1-\frac{13}{3}\,{n}^{-1}+{\frac {
1477}{27}}\,{n}^{-2}-{\frac {93707}{81}}\,{n}^{-3}+{
\frac {8343061}{243}}\,{n}^{-4}-{\frac {2866730137}{
2187}}\,{n}^{-5}+{\frac {1204239422533}{19683}}\,{n}^{
-6} \, + O( n^{-7} ) \right) \quad .
$$
We observe that the ``connective constant'', $8$, is {\it sub-dominant}. With any other initial conditions it would have been
$9$. This is a very rare phenomenon in combinatorics.

The sequence $C_3(n)$ is sequence $A271905$ in the On-Line Encyclopedia of Integer Sequences [Sl].
For the record, here are the first few terms:

1, 6, 37, 240, 1621, 11256, 79717, 572928, 4164841, 30553116, 225817021, 1679454816, 12556853401, 94313192616,  
711189994357, 5381592930816, 40848410792017, 310909645663332, 2372280474687277, 18141232682656320, 
139010366280363601, 1067160872528170536, 8206301850166625797, 63203453697218605440.

We tried to find a recurrence for $C_4(n)$, but, since $160$ terms
did not suffice, we gave up. Nevertheless, using numerics, it
if extremely likely that
$$
C_4(n) \, \sim \, \alpha \, 81^n \cdot n^{-\frac{3}{2}} \quad,
$$
for some constant $\alpha$, but we are unable to conjecture its value.
For the record, here are the first few terms:

1, 24, 997, 51264, 2940841, 180296088, 11559133741, 765337680384,51921457661905, 3590122671128664, 252070718210663749, 17922684123178825536, 
1287832671004683373753, \hfill\break
93368940577497932331288, 6821632357294515590873917, 501741975445243527381995520, \hfill\break
37121266623211130111114816929, 2760712710223967190110979892824, 206267049696409355312012281872181.

The first few terms of $C_5(n)$ are:
1, 120, 44121, 23096640, 14346274601, 9859397817600, 7244702262723241,
5582882474985676800.

The first few terms of $C_6(n)$ are:
1, 720, 2882071, 18754813440, 153480509680141, 1435747717722810960.

Using reliable numeric estimates we are confident in making the following conjecture.

{\bf Conjecture:}
$$
C_d(n) \, \sim \, \alpha_d \, \cdot \, ((d-1)^d)^n \cdot n^{-(d-1)/2} \quad,
$$
for a constant $\alpha_d$. 

One of us (DZ) is pledging \$100 dollars to the OEIS Foundation in honor of the first prover,
and an additional \$25 for an explicit expression for $\alpha_d$ in terms of $d$ and $\pi$.

Readers are welcome to explore further using the Maple package {\tt SVT.txt} available from \hfill\break
{\tt http://www.math.rutgers.edu/\~{}zeilberg/mamarim/mamarimhtml/svt.html}, where there are 
many more terms of the sequences $C_d(n)$ for $3 \leq d \leq 6$.

{\bf References}

[AZ] Moa Apagodu and Doron Zeilberger,
{\it Multi-Variable Zeilberger and Almkvist-Zeilberger Algorithms and the Sharpening of Wilf-Zeilberger Theory},
Adv. Appl. Math. {\bf 37} (2006), 139-152. \hfill\break
{\tt http://www.math.rutgers.edu/\~{}zeilberg/mamarim/mamarimhtml/multiZ.html} \quad .

[FO] Shmuel Friedland and Giorgio Ottaviani, {\it The number of singular vector tuples and uniqueness of best rank-one approximation of tensors},
Journal Found. Comput. Math. {\bf 14} (2014), 1209-1242.  \hfill\break
{\tt http://arxiv.org/abs/1210.8316} \quad .

[M] Major Percy A. MacMahon, {\it Combinatory Analysis}, Originally published by Cambridge University Press, 1917, 1918,
reprinted by Chelsea, N.Y., 1984.

[Sl] Neil J.A. Sloane, ``The On-Line Encyclopedia of Integer Sequences'', {\tt http://www.oeis.org}.

[St] Bernd Sturmfels, {\tt Eigenvectors of Tensors}, Colloquium talk, Rutgers University, April 22, 2016.

[WZ] Jet Wimp and Doron Zeilberger,  {\it Resurrecting the asymptotics of linear recurrences},
J. Math. Anal. Appl. {\bf 111} (1985), 162-177. \hfill\break
{\tt http://www.math.rutgers.edu/\~{}zeilberg/mamarimY/WimpZeilberger1985.pdf} \quad .

[Z] Doron Zeilberger, {\it AsyRec: A Maple package for Computing the Asymptotics of Solutions of Linear Recurrence Equations with Polynomial Coefficients},
The Personal Journal of Shalosh B. Ekhad and Doron Zeilberger, April 4, 2008. \hfill\break
{\tt http://www.math.rutgers.edu/\~{}zeilberg/mamarim/mamarimhtml/asy.html} \quad .

\bigskip
\bigskip
\hrule
\bigskip
Doron Zeilberger, Department of Mathematics, Rutgers University (New Brunswick), Hill Center-Busch Campus, 110 Frelinghuysen
Rd., Piscataway, NJ 08854-8019, USA. \hfill \break
zeilberg at math dot rutgers dot edu \quad ;  \quad {\tt http://www.math.rutgers.edu/\~{}zeilberg/} \quad .
\bigskip
\hrule
\bigskip
Shalosh B. Ekhad, c/o D. Zeilberger, Department of Mathematics, Rutgers University (New Brunswick), Hill Center-Busch Campus, 110 Frelinghuysen
Rd., Piscataway, NJ 08854-8019, USA.
\bigskip
\hrule

\bigskip
Exclusively published in The Personal Journal of Shalosh B. Ekhad and Doron Zeilberger  \hfill \break
({ \tt http://www.math.rutgers.edu/\~{}zeilberg/pj.html})
and {\tt arxiv.org} \quad . 
\bigskip
\hrule
\bigskip
{\bf  First Written:  April 29, 2016} ; 
{\bf  This version:  April 30, 2016} .

\end